\input amstex\documentstyle{amsppt}  
\pagewidth{12.5cm}\pageheight{19cm}\magnification\magstep1
\topmatter
\title On some partitions of a flag manifold\endtitle
\author G. Lusztig\endauthor
\address{Department of Mathematics, M.I.T., Cambridge, MA 02139}\endaddress
\thanks{Supported in part by the National Science Foundation}\endthanks
\endtopmatter   
\document

\define\uWW{\un\WW}

\define\hT{\hat T}

\define\Lie{\text{\rm Lie }}

\define\sqc{\sqcup}

\define\bZ{\bar Z}

\define\op{\oplus}

\define\part{\partial}

\define\m{\mapsto}
\define\do{\dots}

\define\lra{\leftrightarrow}

\define\sub{\subset}    

\define\T{\times}
\define\ti{\tilde}
\define\nl{\newline}
\redefine\i{^{-1}}

\define\un{\underline}
\define\ov{\overline}
\define\ot{\otimes}
\define\bbq{\bar{\QQ}_l}

\define\tr{\text{\rm tr}}

\define\g{\gamma}

\define\e{\epsilon}

\define\ph{\phi}

\define\r{\rho}
\define\s{\sigma}

\define\th{\theta}

\define\x{\xi}

\define\Ph{\Phi}

\define\kk{\bold k}

\redefine\ss{\bold s}

\define\CC{\bold C}

\define\FF{\bold F}

\define\NN{\bold N}

\define\QQ{\bold Q}

\define\SS{\bold S}

\define\WW{\bold W}
\define\ZZ{\bold Z}

\define\cb{\Cal B}
\define\cc{\Cal C}

\define\cf{\Cal F}

\define\ch{\Cal H}

\define\co{\Cal O}

\define\ct{\Cal T}

\define\cw{\Cal W}

\define\fg{\frak g}

\define\fk{\frak k}

\define\fp{\frak p}

\define\Ir{\text{\rm Irr}}

\define\CR{CR}
\define\DL{DL}
\define\GP{GP}
\define\KL{KL}
\define\RE{L1}
\define\UE{L2}
\define\CS{L3}
\define\CV{L4}
\define\AS{L5}
\define\SH{Sh}
\define\SP{S1}
\define\SPP{S2}

\subhead 1\endsubhead
Let $G$ be a connected reductive group over an algebraically closed field $\kk$ of characteristic exponent 
$p\ge1$. We shall assume that $p$ is $1$ or a good prime for $G$. Let $\cb$ be the variety of Borel subgroups of 
$G$ and let $\WW$ be the Weyl group of $G$. Note that $\WW$ naturally indexes ($w\m\co_w$) the orbits of $G$ 
acting on $\cb\T\cb$ by simultaneous conjugation on the two factors. For $g\in G$ we set 
$\cb_g=\{B\in\cb;g\in B\}$. The varieties $\cb_g$ play an important role in representation theory and their 
geometry has been studied extensively. More generally for $g\in G$ and $w\in\WW$ we set 
$$\cb_g^w=\{B\in\cb;(B,gBg\i)\in\co_w\}.$$
Note that $\cb_g^1=\cb_g$ and that for fixed $g$, $(\cb_g^w)_{w\in\WW}$ form a partition of the flag manifold 
$\cb$.

For fixed $w$, the varieties $\cb_g^w$ ($g\in G$) appear as fibres of a map to $G$ which was introduced in 
\cite{\CS} as part of
the definition of character sheaves. Earlier, the varieties $\cb_g^w$ for $g$ regular semisimple appeared in 
\cite{\RE} (a precursor of \cite{\CS}) where it was shown that from their topology (for $\kk=\CC$) one can extract
nontrivial information about the character table of the corresponding group over a finite field. 

In this paper we present some experimental evidence for our belief that from the topology of $\cb_g^w$ for $g$ 
unipotent (with $\kk=\CC$) one can extract information closely connected with the map \cite{\KL, 9.1} from 
unipotent classes in $G$ to conjugacy classes in $\WW$. 

I thank David Vogan for some useful discussions.

\subhead 2\endsubhead
We fix a prime number $l$ invertible in $\kk$. Let $g\in G$ and $w\in\WW$. For $i,j\in\ZZ$ let 
$H^i_c(\cb_g^w,\bbq)_j$ be the subquotient of pure weight $j$ of the $l$-adic cohomology space 
$H^i_c(\cb_g^w,\bbq)$. The centralizer $Z(g)$ of $g$ in $G$ acts on $\cb_g^w$ by conjugation and this induces an 
action of the group of components $\bZ(g)$ on $H^i_c(\cb_g^w,\bbq)$ and on each $H^i_c(\cb_g^w,\bbq)_j$. For 
$z\in\bZ(g)$ we set
$$\Xi^w_{g,z}=\sum_{i,j\in\ZZ}(-1)^i\tr(z,H^i_c(\cb_g^w,\bbq)_j)v^j\in\ZZ[v]$$
where $v$ is an indeterminate; the fact that this belongs to $\ZZ[v]$ and is independent of the choice of $l$ is 
proved by an argument similar to that in the proof of \cite{\DL, 3.3}. 

Let $l:\WW@>>>\NN$ be the standard length function. The simple reflections $s\in\WW$ (that is the elements of
length $1$ of $\WW$) are numbered as $s_1,s_2,\do$. Let $w_0$ be the element of maximal length in $\WW$. Let 
$\uWW$ be the set of conjugacy classes in $\WW$.

Let $\ch$ be the Iwahori-Hecke algebra of $\WW$ with parameter $v^2$ (see \cite{\GP, 4.4.1}; in the definition in 
{\it loc.cit.} we take $A=\ZZ[v,v^{-1}],a_s=b_s=v^2$). Let $(T_w)_{w\in\WW}$ be the standard basis of $\ch$ (see
\cite{\GP, 4.4.3, 4.4.6}). For $w\in\WW$ let $\hT=v^{-2l(w)}T_w$. If $s_{i_1}s_{i_2}\do s_{i_t}$ is a reduced 
expression for $w\in\WW$ we write also $\hT_w=\hT_{i_1i_2\do i_t}$.

For any $g\in G,z\in\bZ(g)$ we set 
$$\Pi_{g,z}=\sum_{w\in\WW}\Xi^w_{g,z}\hT_w\in\ch.$$
The following result can be proved along the lines of the proof of \cite{\DL, Theorem 1.6} (we replace the
Frobenius map in that proof by conjugation by $g$); alternatively, for $g$ unipotent, we may use 6(a).

(a) {\it $\Pi_{g,z}$ belongs to the centre of the algebra $\ch$.}
\nl
According to \cite{\GP, 8.2.6, 7.1.7}, an element $c=\sum_{w\in\WW}c_w\hT_w$ ($c_w\in\ZZ[v,v\i]$) in the centre of
$\ch$ is uniquely determined by the coefficients $c_w (w\in\WW_{min})$ and we have $c_w=c_{w'}$ if
$w,w'\in\WW_{min}$ are conjugate in $\WW$; here $\WW_{min}$ is the set of elements of $\WW$ which have minimal 
length in their conjugacy class. This applies in particular to $c=\Pi_{g,z}$, see (a). For any $C\in\uWW$ we set 
$\Xi^C_{g,z}=\Xi^w_{g,z}$ where $w$ is any element of $C\cap\WW_{min}$.

Note that if $g=1$ then $\Pi_{g,1}=(\sum_wv^{2l(w)}1$. If $g$ is regular unipotent then 
$\Pi_{g,1}=\sum_{w\in\WW}v^{2l(w)}\hT_w$. If $G=PGL_3(\kk)$ and $g\in G$ is regular semisimple then
$\Pi_{g,1}=6+3(v^2-1)(\hT_1+\hT_2)+(v^2-1)^2(\hT_{12}+\hT_{21})+(v^6-1)\hT_{121}$; if $g\in G$ is a transvection
then $\Pi_{g,1}=(2v^2+1)+v^4(\hT_1+\hT_2)+v^6\hT_{121}$.

For $g\in G$ let $cl(g)$ be the $G$-conjugacy class of $g$; let $\ov{cl(g)}$ be the closure of $cl(g)$. Let 
$\SS_g$ be the set of all $C\in\uWW$ such that $\Xi^C_{g,1}\ne0$ and $\Xi^C_{g',1}=0$ for any 
$g'\in\ov{cl(g)}-cl(g)$. If $\cc$ is a conjugacy class in $G$ we shall also write $\SS_\cc$ instead of $\SS_g$ 
where $g\in\cc$.

We describe the set $\SS_g$ and the values $\Xi^C_{g,1}$ for $C\in\SS_g$ for various $G$ of low rank and various
unipotent elements $g$ in $G$. We denote by $u_n$ a unipotent element of $G$ such that $\dim\cb_{u_n}=n$. The 
conjugacy class of $w\in\WW$ is denoted by $(w)$.

$G$ of type $A_1$.
$$\SS_{u_1}=(1), \SS_{u_0}=(s_1);\Xi^1_{u_1,1}=1+v^2, \Xi^{s_1}_{u_0,1}=v^2.$$
$G$ of type $A_2$.
$$\SS_{u_3}=(1), \SS_{u_1}=(s_1), \SS_{u_0}=(s_1s_2).$$
$$\Xi^1_{u_3,1}=1+2v^2+2v^4+v^6, \Xi^{(s_1)}_{u_1,1}=v^4,\Xi^{(s_1s_2)}_{u_0,1}=v^4.$$
$G$ of type $B_2$. (The simple reflection corresponding to the long root is denoted by $s_1$.)
$$\SS_{u_4}=(1), \SS_{u_2}=(s_1), \SS_{u_1}=\{(s_2),(s_1s_2s_1s_2)\}, \SS_{u_0}=(s_1s_2).$$
$$\align&\Xi^1_{u_4,1}=(1+v^2)^2(1+v^4), \Xi^{(s_1)}_{u_2,1}=v^4(1+v^2), \Xi^{(s_2)}_{u_1,1}=2v^4,\\&
\Xi^{(s_1s_2s_1s_2)}_{u_1,1}=v^6(v^2-1), \Xi^{(s_1s_2)}_{u_0,1}=v^4.\endalign$$ 
$G$ of type $G_2$. (The simple reflection corresponding to the long root is denoted by $s_2$.)
$$\SS_{u_6}=(1), \SS_{u_3}=(s_2), \SS_{u_2}=\{(s_1),(s_1s_2s_1s_2s_1s_2)\}, \SS_{u_1}=(s_1s_2s_1s_2),
\SS_{u_0}=(s_1s_2).$$
$$\align&\Xi^1_{u_6,1}=(1+v^2)^2(1+v^4+v^8), \Xi^{(s_2)}_{u_3,1}=v^6(1+v^2), \Xi^{(s_1)}_{u_2,1}=v^4(1+v^2),\\&
\Xi^{(s_1s_2s_1s_2s_1s_2)}_{u_2,1}=v^8(v^4-1), \Xi^{(s_1s_2s_1s_2)}_{u_1,1}=2v^8, \Xi^{(s_1s_2)}_{u_0,1}=v^4.
\endalign$$ 
$G$ of type $B_3$. (The simple reflection corresponding to the short root is denoted by $s_3$ and $(s_1s_3)^2=1$.)
$$\align&\SS_{u_9}=(1), \SS_{u_5}=(s_1), \SS_{u_4}=\{(s_3),(s_2s_3s_2s_3)\}, \SS_{u_3}=\{(s_1s_3),(w_0)\},\\&
\SS_{u_2}=(s_1s_2), \SS_{u_1}=\{(s_2s_3),(s_2s_3s_1s_2s_3)\}, \SS_{u_0}=(s_1s_2s_3).\endalign$$
$$\align&\Xi^1_{u_9,1}=(1+v^2)^3(1+v^4)(1+v^4+v^8), \Xi^{(s_1)}_{u_5,1}=v^8(1+v^2)^2, \\&
\Xi^{(s_2s_3s_2s_3)}_{u_4,1}=v^8(1+v^2)(v^4-1), \Xi^{(s_3)}_{u_4,1}=2v^6(1+v^2)^2,\\&
\Xi^{(s_1s_3)}_{u_3,1}=v^8(1+v^2), \Xi^{(w_0)}_{u_3,1}=v^{14}(v^4-1),\Xi^{(s_1s_2)}_{u_2,1}=2v^8, \\&
\Xi^{(s_2s_3)}_{u_1,1}=2v^6,\Xi^{(s_2s_3s_1s_2s_3)}_{u_1,1}=v^8(v^2-1), \Xi^{(s_1s_2s_3)}_{u_0,1}=v^6.\endalign$$
$G$ of type $C_3$. (The simple reflection corresponding to the long root is denoted by $s_3$ and $(s_1s_3)^2=1$;
$u''_2$ denotes a unipotent element which is regular inside a Levi subgroup of type $C_2$; $u'_2$ denotes
a unipotent element with $\dim\cb_{u''_2}=2$ which is not conjugate to $u''_2$.) 
$$\align&\SS_{u_9}=(1), \SS_{u_6}=(s_3), \SS_{u_4}=\{(s_1),(s_2s_3s_2s_3)\}, \SS_{u_3}=\{(s_1s_3),(w_0)\},\\&
\SS_{u'_2}=(s_1s_2), \SS_{u''_2}=(s_2s_3), \SS_{u_1}=(s_2s_3s_1s_2s_3), \SS_{u_0}=(s_1s_2s_3).\endalign$$
$$\align&\Xi^1_{u_9,1}=(1+v^2)^3(1+v^4)(1+v^4+v^8), \Xi^{(s_3)}_{u_6,1}=v^6(1+v^2)^2(1+v^4),\\&
\Xi^{(s_2s_3s_2s_3)}_{u_4,1}=v^{10}(v^4-1), \Xi^{(s_1)}_{u_4,1}=2v^8(1+v^2),\\&
\Xi^{(s_1s_3)}_{u_3,1}=v^8(1+v^2), \Xi^{(w_0)}_{u_3,1}=v^{14}(v^4-1),\Xi^{(s_1s_2)}_{u'_2,1}=v^6(1+v^2),\\& 
\Xi^{(s_2s_3)}_{u''_2,1}=v^6(1+v^2),
\Xi^{(s_2s_3s_1s_2s_3)}_{u_1,1}=v^{10}, \Xi^{(s_1s_2s_3)}_{u_0,1}=v^6.\endalign$$

\subhead 3\endsubhead
Consider the following properties of $G$:
$$\WW=\sqc_u\SS_u\tag a$$
($u$ runs over a set of representatives for the unipotent classes in $G$); 
$$c_u\in\SS_u\tag b$$
where for any unipotent element $u\in G$, $c_u$ denotes the conjugacy class in $\WW$  associated to $u$ in 
\cite{\KL, 9.1}. 
Note that (a),(b) hold for $G$ of rank $\le3$ (see \S2). We have used the computations of
the map in \cite{\KL, 9.1} given in \cite{\KL, \S9}, \cite{\SP}, \cite{\SPP}. 
We will show elsewhere that (a), (b) hold for $G$ of type $A_n$ or $C_n$. We expect that a proof similar
to that in type $C_n$ would yield (a),(b) in type $B_n$ and $D_n$.
We expect that (a),(b) hold also for $G$ simple of exceptional type; (a) should follow by computing
the product of some known (large) matrices using 6(a).
The equality $\WW=\cup_u\SS_u$ is clear since for a regular unipotent $u$ and any $w$ we have 
$\Xi_{u,1}^w=v^{2l(w)}$.  

\subhead 4\endsubhead
Assume that $G=Sp_{2n}(\kk)$. The Weyl group $\WW$ can be identified in the standard way with the subgroup of the
symmetric group $S_{2n}$ consisting of all permutations of $\{1,2,\do,n,n',\do,2',1'\}$ which commute with the 
involution $1\lra1', 2\lra2',\do, n\lra n'$. We say that two elements of $\uWW$ are equivalent if they are
contained in the same conjugacy class of $S_{2n}$. The set of equivalence classes in $\uWW$ is in bijection with 
the set of partitions of $2n$ in which every odd part appears an even number of times (to $C\in\uWW$ we attach the
partition which has a part $j$ for every $j$-cycle of an element of $C$ viewed as a permutation of 
$\{1,2,\do,n,n',\do,2',1'\}$). The same set of partitions of $2n$ indexes the set of unipotent classes of $G$.
Thus we obtain a bijection between the set of equivalence classes in $\uWW$ and the set of unipotent classes of
$G$. In other words we obtain a surjective map $\ph$ from $\uWW$ to the set of unipotent classes of $G$ whose
fibres are the equivalence classes in $\uWW$. One can show that
for any unipotent class $\cc$ in $G$ we have $\ph\i(\cc)=\SS_u$ where $u\in\cc$. 

\subhead 5\endsubhead
Recall that the set of unipotent elements in $G$ can be partitioned into "special pieces" (see \cite{\AS}) where
each special piece is a union of unipotent classes exactly one of which is "special". Thus the special pieces can
be indexed by the set of isomorphism classes of special representations of $\WW$ which depends only on $\WW$ as a
Coxeter group (not on the underlying root system). For each special piece $\s$ of $G$ we consider the subset 
$\SS_\s:=\sqc_{\cc\sub\s}\SS_\cc$ of $\uWW$ (here $\cc$ runs over the unipotent classes contained in $\s$). We
expect that each such subset $\SS_\s$ depends only on the Coxeter group structure of $\WW$ (not on the underlying
root system). As evidence for this we note that the subsets $\SS_\s$ for $G$ of type $B_3$ are the same as the 
subsets $\SS_\s$ for $G$ of type $C_3$. These subsets are as follows: 
$$\align&\{1\}, \{(s_1),(s_3),(s_2s_3s_2s_3)\}, \{(s_1s_3),(w_0)\}, \{(s_1s_2)\},\\&
\{(s_2s_3),(s_2s_3s_1s_2s_3)\}, \{(s_1s_2s_3)\}.\endalign$$

\subhead 6\endsubhead
Let $g\in G$ be a unipotent element and let $z\in\bZ(g)$, $w\in W$. We show how the polynomial $\Xi_{g,z}^w$ can 
be computed using information from representation theory. We may assume that $p>1$ and that $\kk$ is the algebraic
closure of the finite field $\FF_p$. We choose an $\FF_p$ split rational structure on $G$ with Frobenius map 
$F_0:G@>>>G$. We may assume that $g\in G^{F_0}$. Let $q=p^m$ where $m\ge1$ is sufficiently divisible. In 
particular $F:=F_0^m$ acts trivially on $\bZ(g)$ hence $cl(g)^F$ is a union of $G^F$-conjugacy classes naturally 
indexed by the conjugacy classes in $\bZ(g)$; in particular the $G^F$-conjugacy class of $g$ corresponds to 
$1\in\bZ(g)$. Let $g_z$ be an element of the $G^F$-conjugacy class in $cl(g)^F$ corresponding to the 
$\bZ(g)$-conjugacy class of $z\in\bZ(g)$. The set $\cb_{g_z}^w$ is $F$-stable. We first compute the number of 
fixed points $|(\cb_{g_z}^w)^F|$. 

Let $\ch_q=\bbq\ot_{\ZZ[v,v\i]}\ch$ where $\bbq$ is regarded as a $\ZZ[v,v\i]$-algebra with $v$ acting as 
multiplication by $\sqrt{q}$. We write $T_w$ instead of $1\ot T_w$. Let $\Ir\WW$ be a set of representatives for
the isomorphism classes of irreducible $\WW$-modules over $\bbq$. For any $E\in\Ir\WW$ let $E_q$ be the 
irreducible $\ch_q$-module corresponding naturally to $E$. Let $\cf$ be the vector space of functions 
$\cb^F@>>>\bbq$. We regard $\cf$ as a $G^F$-module by $\g:f\m f'$, $f'(B)=f(\g\i B\g)$ for all $B\in\cb^F$. We 
identify $\ch_q$ with the algebra of all endomorphisms of $\cf$ which commute with the $G^F$-action, by 
identifying $T_w$ with the endomorphism $f\m f'$ where $f'(B)=\sum_{B'\in\cb^F;(B,B')\in\co_w}f(B)$ for all 
$B\in\cb^F$. As a module over $\bbq[G^F]\ot\ch_q$ we have canonically $\cf=\op_{E\in\Ir\WW}\r_E\ot E_q$ where 
$\r_E$ is an irreducible $G^F$-module. Hence if $\g\in G^F$ and $w\in\WW$ we have 
$\tr(\g T_w,\cf)=\sum_{E\in\Ir\WW}\tr(\g,\r_E)\tr(T_w,E_q)$. From the definition we have
$\tr(\g T_w,\cf)=|\{B\in\cb^F;(B,\g B\g\i)\in\co_w\}|=|(\cb_\g^w)^F|$. Taking $\g=g_z$ we obtain
$$|(\cb_{g_z}^w)^F|=\sum_{E\in\Ir\WW}\tr(g_z,\r_E)\tr(T_w,E_q).\tag a$$ 
The quantity $\tr(g_z,\r_E)$ can be computed explicitly, by the method of \cite{\CV}, in terms of generalized 
Green functions and of the entries of the non-abelian Fourier transform matrices \cite{\UE}; in particular it is a
polynomial with rational coefficients in $\sqrt{q}$. The quantity $\tr(T_w,E_q)$ can be also computed explicitly
(see \cite{\GP}, Ch.10,11); it is a polynomial with integer coefficients in $\sqrt{q}$. Thus $|(\cb_{g_z}^w)^F|$ 
is an explicitly computable polynomial with rational coefficients in $\sqrt{q}$. Substituting here $\sqrt{q}$ by 
$v$ we obtain the polynomial $\Xi_{g,z}^w$. This argument shows also that $\Xi_{g,z}^w$ is independent of $p$ 
(note that the pairs $(g,z)$ up to conjugacy may be parametrized by a set independent of $p$).

This is how the various $\Xi_{g,z}^w$ in \S2 were computed, except in type $A_1,A_2,B_2$ where they were computed 
directly from the definitions. (For type $B_3,C_3$ we have used the computation of Green functions in \cite{\SH};
for type $G_2$ we have used directly \cite{CR} for the character of $\r_E$ at unipotent elements.)

\subhead 7\endsubhead
In this section we assume that $G$ is simply connected. Let $\ti G=G(\kk((\e)))$ where $\e$ is an indeterminate.
Let $\ti\cb$ be the set of Iwahori subgroups of $\ti G$. Let $\ti\WW$ the affine Weyl group attached to $\ti G$. 
Note that $\ti\WW$ naturally indexes ($w\m\co_w$) the orbits of $\ti G$ acting on $\ti\cb\T\ti\cb$ by simultaneous
conjugation on the two factors. For $g\in\ti G$ and $w\in\ti\WW$ we set 
$$\ti\cb_g^w=\{B\in\ti\cb;(B,gBg\i)\in\co_w\}.$$
By analogy with \cite{\KL, \S3} we expect that when $g$ is regular semisimple, $\ti\cb_g^w$ has a natural 
structure of a locally finite union of algebraic varieties over $\kk$ of bounded dimension and that, moreover, if
$g$ is also elliptic, then $\ti\cb_g^w$ has a natural structure of algebraic variety over $\kk$.
It would follow that for $g$ elliptic and $w\in\ti\WW$,
$$\Xi^w_g=\sum_{i,j\in\ZZ}(-1)^i\dim H^i_c(\ti\cb_g^w,\bbq)_jv^j\in\ZZ[v]$$
is well defined; one can then show that the formal sum $\sum_{w\in\ti\WW}\Xi^w_g\hT_w$
is central in the completion of the affine Hecke algebra consisting of all formal sums 
$\sum_{w\in\ti\WW}a_w\hT_w$ ($a_w\in\QQ(v)$) that is, it commutes with any $\hT_w$. (Here $\hT_w$ is defined
as in \S2 and the completion of the affine Hecke algebra is regarded as a bimodule over the actual affine Hecke
algebra in the natural way.)

\subhead 8\endsubhead
In this and the next subsection we assume that $G$ is adjoint. For $g\in G,z\in\bZ(g),w\in\WW$ we set
$$\xi^w_{g,z}=\Xi^w_{g,z}|_{v=1}=\sum_{i\in\ZZ}(-1)^i\tr(z,H^i_c(\cb_g^w,\bbq))\in\ZZ.$$
This integer is independent of $l$. For any $g\in G,z\in\bZ(g)$ we set 
$$\pi_{g,z}=\sum_{w\in\WW}\xi^w_{g,z}w\in\ZZ[W].$$
This is the specialization of $\Pi_{g,z}$ for $v=1$. Hence from 2(a) we see that $\pi_{g,z}$ is in the centre of 
the ring $\ZZ[\WW]$. Thus for any $C\in\uWW$ we can set $\xi^C_{g,z}=\xi^w_{g,z}$ where $w$ is any element of
$C$. For $g\in G$ let $\ss_g$ be the set of all $C\in\uWW$ such that $\xi^C_{g,z}\ne0$ for some $z\in\bZ(g)$ and 
$\xi^C_{g',z'}=0$ for any $g'\in\ov{cl(g)}-cl(g)$ and any $z'\in\bZ(g')$. We describe the set $\ss_g$ and the
values $\xi^C_{g,z}=0$ for $C\in\ss_g$, $z\in\bZ(g)$, for various $G$ of low rank and various unipotent elements 
$g$ in $G$. We use the notation in \S2. Moreover in the case where $\bZ(g)\ne\{1\}$ we denote by $z_n$ an element 
of order $n$ in $\bZ(g)$.

$G$ of type $A_1$.
$$\ss_{u_1}=(1), \ss_{u_0}=(s_1);\xi^1_{u_1,1}=2, \xi^{s_1}_{u_0,1}=1.$$
$G$ of type $A_2$.
$$\ss_{u_3}=(1), \ss_{u_1}=(s_1), \ss_{u_0}=(s_1s_2).$$
$$\xi^1_{u_3,1}=6, \xi^{(s_1)}_{u_1,1}=1,\xi^{(s_1s_2)}_{u_0,1}=1.$$
$G$ of type $B_2$. 
$$\ss_{u_4}=(1), \ss_{u_2}=(s_1), \ss_{u_1}=\{(s_2),(s_1s_2s_1s_2)\}, \ss_{u_0}=(s_1s_2).$$
$$\align&\xi^1_{u_4,1}=8, \xi^{(s_1)}_{u_2,1}=2, \xi^{(s_2)}_{u_1,1}=2, \xi^{(s_1s_2s_1s_2)}_{u_1,1}=0, \\&
\xi^{(s_2)}_{u_1,z_2}=0, \xi^{(s_1s_2s_1s_2)}_{u_1,z_2}=2, \xi^{(s_1s_2)}_{u_0,1}=1.\endalign$$ 
$G$ of type $G_2$. 
$$\ss_{u_6}=(1), \ss_{u_3}=(s_2), \ss_{u_2}=(s_1), \ss_{u_1}=\{(s_1s_2s_1s_2s_1s_2),(s_1s_2s_1s_2)\},
\ss_{u_0}=(s_1s_2).$$
$$\align&\xi^1_{u_6,1}=12, \xi^{(s_2)}_{u_3,1}=2, \xi^{(s_1)}_{u_2,1}=2,
\xi^{(s_1s_2s_1s_2s_1s_2)}_{u_1,1}=-3, \xi^{(s_1s_2s_1s_2s_1s_2)}_{u_1,z_2}=3, \\&
\xi^{(s_1s_2s_1s_2s_1s_2)}_{u_1,z_3}=0, \xi^{(s_1s_2s_1s_2)}_{u_1,1}=2, 
\xi^{(s_1s_2s_1s_2)}_{u_1,z_2}=0, \xi^{(s_1s_2s_1s_2)}_{u_1,z_3}=2, \xi^{(s_1s_2)}_{u_0,1}=1.\endalign$$
$G$ of type $B_3$. 
$$\align&\ss_{u_9}=(1), \ss_{u_5}=(s_1), \ss_{u_4}=\{(s_3),(s_2s_3s_2s_3)\}, \ss_{u_3}=(s_1s_3),\\&
\ss_{u_2}=\{(s_1s_2),(w_0)\}, \ss_{u_1}=\{(s_2s_3),(s_2s_3s_1s_2s_3)\}, \ss_{u_0}=(s_1s_2s_3).\endalign$$
$$\align&\xi^1_{u_9,1}=48, \xi^{(s_1)}_{u_5,1}=4, \xi^{(s_2s_3s_2s_3)}_{u_4,1}=0,
 \xi^{(s_2s_3s_2s_3)}_{u_4,z_2}=4, \xi^{(s_3)}_{u_4,1}=8,\\&
 \xi^{(s_3)}_{u_4,1}=0,\xi^{(s_1s_3)}_{u_3,1}=2,
 \xi^{(w_0)}_{u_2,1}=0, \xi^{(w_0)}_{u_2,z_2}=6\\&
\xi^{(s_1s_2)}_{u_2,1}=2, \xi^{(s_1s_2)}_{u_2,z_2}=0, 
\xi^{(s_2s_3)}_{u_1,1}=2,\xi^{(s_2s_3)}_{u_1,z_2}=0,\\&
\xi^{(s_2s_3s_1s_2s_3)}_{u_1,1}=0, \xi^{(s_2s_3s_1s_2s_3)}_{u_1,z_2}=2, \xi^{(s_1s_2s_3)}_{u_0,1}=1.\endalign$$
$G$ of type $C_3$. 
$$\align&\ss_{u_9}=(1), \ss_{u_6}=(s_3), \ss_{u_4}=\{(s_1),(s_2s_3s_2s_3)\}, \ss_{u_3}=(s_1s_3),\\&
\ss_{u'_2}=(s_1s_2), \ss_{u''_2}=(s_2s_3), \ss_{u_1}=\{(s_2s_3s_1s_2s_3),w_0\} \ss_{u_0})=(s_1s_2s_3).\endalign$$
$$\align&\xi^1_{u_9,1}=48, \xi^{(s_3)}_{u_6,1}=8,
\xi^{(s_2s_3s_2s_3)}_{u_4,1}=0, \xi^{(s_2s_3s_2s_3)}_{u_4,z_2}=4, \\&
\xi^{(s_1)}_{u_4,1}=4, \xi^{(s_1)}_{u_4,1}=0,\xi^{(s_1s_3)}_{u_3,1}=2, 
\xi^{(s_1s_2)}_{u'_2,1}=2, \xi^{(s_2s_3)}_{u''_2,1}=2,\\&
\xi^{(s_2s_3s_1s_2s_3)}_{u_1,1}=1, \xi^{(s_2s_3s_1s_2s_3)}_{u_1,z_2}=1, 
\xi^{(w_0)}_{u_1,1}=-3, \xi^{(w_0)}_{u_1,z_2}=3, \xi^{(s_1s_2s_3)}_{u_0,1}=1.\endalign$$

\subhead 9\endsubhead
For any unipotent element $u\in G$ let $n_u$ be the number of isomorphism classes of irreducible representations 
of $\bZ(u)$ which appear in the Springer correspondence for $G$. Consider the following properties of $G$:
$$\WW=\sqc_u\ss_u\tag a$$
($u$ runs over a set of representatives for the unipotent classes in $G$); for any unipotent element $u\in G$,
$$|\ss_u|=n_u.\tag b$$
By the results in \S8, (a,(b) hold if $G$ has rank $\le3$. We will show elsewhere that (a),(b) hold if $G$ is of 
type $A$. We expect that (a),(b) hold in general. 
The equality $\WW=\cup_u\ss_u$ is clear since for a regular unipotent $u$ and any $w$ we have $\xi_{u,1}^w=1$. 
Consider also the following property of $G$: for any $g\in G$, $w\in\WW$,
$$\align&\text{$\xi_{g,1}^w$ is equal to the trace of $w$ on the Springer representation} \\&
\text{of $\WW$ on $\op_iH^{2i}(\cb_g,\bbq)$.}\tag c\endalign$$
Again (c) holds if $G$ is of type $A$ and in the examples in \S8; we expect
that it holds in general. Note that in (c) one can ask whether for any $z$,
$\xi_{g,z}^w$ is equal to the trace of $wz$ on the Springer representation of $\WW\T\bZ(g)$ on 
$\op_iH^{2i}(\cb_g,\bbq)$; but such an equality is not true in general for $z\ne1$ (for example for $G$ of type
$B_2$). 

\subhead 10\endsubhead
In this subsection we assume that $\kk=\CC$. We show how the method of \cite{\KL, 9.1} extends to the case of
symmetric spaces.

Let $K$ be a subgroup of $G$ which is the fixed point set of an involution $\th:G@>>>G$.
Let $\fg,\fk$ be the Lie algebras of $G,K$. Let $\fp$ be the $(-1)$-eigenspace of $\th:\fg@>>>\fg$.
Let $\fp_{nil}$ be the set of elements $x\in\fp$ which are nilpotent in $\fg$. Note that $K$ acts naturally
on $\fp_{nil}$ with finitely many orbits.
Let $\ct$ be the variety of all tori $T$ in $G$ such that $\th(t)=t\i$ for all $t\in T$ and such that $T$
has maximum possible dimension. It is known that $K$ acts transitively on $\ct$ by conjugation.
For $T\in\ct$ let $\cw_T$ be the normalizer of $T$ in $K$ modulo the centralizer 
of $T$ in $K$. Let $\un{\cw}$ be the set of conjugay classes in the finite group $\cw_T$; this is independent
of the choice of $T$.

Let $\Ph=\CC((\e))$, $A=\CC[[\e]]$ where $\e$ is an indeterminate. Let $\Ph'$ be an algebraic closure of $\Ph$.
Then the groups $G(\Ph),K(\Ph)$ are well defined and the set $\ct(\Ph)$ of $\Ph$-points of $\ct$ is well
defined. Let $\fp_\Ph=\Ph\ot\fp$, $\fp_A=A\ot\fp$.

The group $K(\Ph)$ acts naturally by conjugation on $\ct(\Ph)$; as in \cite{\KL, \S1, Lemma 2} we see that
the set of $K(\Ph)$-orbits on $\ct(\Ph)$ is naturally in $1-1$ correspondence with the set $\un{\cw}$.
For $\g\in\un{\cw}$ let $\co_\g$ be the $K(\Ph)$-orbit on $\ct(\Ph)$ corresponding to $\g$.

An element $\x\in\fp_\Ph$ is said to be "regular semisimple" if there is a unique $T'\in\ct(\Ph')$ such that
$\x\in\Lie(T')$; we then set $T'_\x=T'$ and we have necessarily $T'_\x\in\ct(\Ph)$.

Let $N\in\fp_{nil}$. We consider the subset $N+\e\fp_A$ of $\fp_\Ph$.
As in \cite{\KL, 9.1} there is a unique element $\g\in\un{\cw}$ such that the following holds: 
there exists a "Zariski open dense" subset $V$ of $N+\e\fp_A$ such that
for any $\x\in V$, $\x$ is "regular semisimple" and $T'_\x\in\co_\g$.
Note that $N\m\g$ is constant on $K$-orbits hence it defines a map $\Psi$ from the set of
$K$-orbits on $\fp_{nil}$ to $\un{\cw}$.
(In the case where $G=H\T H$ where $H$ is a connected reductive group, $K$ is the diagonal in $H\T H$ and 
$\th(a,b)=(b,a)$, the map $\Psi$ reduces to the map defined in \cite{\KL, 9.1}.)

One can show that if $(G,K)=(GL_{2n}(\CC),Sp_{2n}(\CC))$, then $\Psi$ is a bijection.
Note that for general $(G,K)$, $\Psi$ is neither injective nor surjective.

\enddocument